\documentstyle[12pt]{article}
\def\Bbb R{{\rm \bf R}}
\def\proclaim#1{\vskip2mm{\bf #1}\em}
\def\endproclaim{\em \vskip2mm}
\def\tag#1{\eqno(#1)}
\def\gathered{\begin{array}{c}}
\def\endgathered{\end{array}}
\def\text{\mbox}

\begin{document}

\title {A remark on the enclosure method for a body with an unknown homogeneous background conductivity}
\author{Masaru IKEHATA\\
Department of Mathematics,
Graduate School of Engineering\\
Gunma University, Kiryu 376-8515, JAPAN}
%\date{11 October 2007}
\maketitle

\begin{abstract}
Previous applications of the enclosure method with a {\it finite} set of observation data to a mathematical
model of electrical impedance tomography are based on the
assumption that the conductivity of the background body is homogeneous and
{\it known}. This paper considers the case when the conductivity
is homogeneous and {\it unknown}. It is shown that,
in two dimensions if the domain
occupied by the background body is enclosed by an {\it ellipse},
then it is still possible to extract some information about the location
of unknown cavities or inclusions embedded in the body without
knowing the background conductivity provided the Fourier series
expansion of the voltage on the boundary does not
contain high frequency parts (band limited) and satisfies a non
vanishing condition of a quantity involving the Fourier
coefficients.

\noindent
AMS: 35R30

\end{abstract}

\section{Introduction}

The aim of this paper is to reconsider previous applications
\cite{I1,Ie2} of the {\it enclosure method} with a {\it finite}
set of observation data to inverse boundary value problems related
to a continuum model of electrical impedance tomography
\cite{B1,B2}. The point is: those applications are based on the
assumption that the conductivity of the background body is
homogeneous and {\it known}. However, from a mathematical point of
view, the problem whether or not one can still extract some
information about unknown discontinuity from the finite set of
observation data without knowing the exact value of the
conductivity is quite interesting.  Proofs of some previous known
uniqueness results that employ a finite set of observation data,
for example, \cite{FV} for cracks and \cite{FI, S} for inclusions
are based on the assumption that the conductivity of the
background body is known.  This is because they start with
applying the uniqueness of the Cauchy problem for elliptic
equations.

Besides needless to say, we cannot know the exact value of the conductivity of the background body.
The inaccurate value causes an error on the observation data and therefore on the indicator
function in the enclosure method.

In order to describe the problem more precisely let us start with
recalling a typical application of the enclosure method with a
single set of observation data.

\noindent
Let $\Omega$ be a bounded domain of $\Bbb R^2$ with Lipschitz boundary.
Let $D$ be an open subset with Lipschitz boundary of $\Omega$ such that
$\overline D\subset\Omega$ and $\Omega\setminus\overline D$ is connected.
Consider a non constant solution of the elliptic problem:
$$\begin{array}{c}
\displaystyle
\triangle u=0\,\,\text{in}\,\Omega\setminus\overline D,\\
\\
\displaystyle
\frac{\partial u}{\partial\nu}=0\,\,\text{on}\,\partial D.
\end{array}
\tag {1.1}
$$
\noindent
Here $\nu=(\nu_1,\nu_2)$ denotes the unit outward normal vector field on $\partial(\Omega\setminus\overline D)$.
The $D$ is a mathematical model of the union of cavities inside the body.

In \cite{I1} we considered the problem of extracting information
about the location and shape of $D$ in two dimensions from the
observation data that is a single set of Cauchy data of $u$ on
$\partial\Omega$.   Assuming that $D$ is given by the inside of a
polygon with an additional condition on the diameter, we
established an extraction formula of the {\it convex hull} of $D$
from the data.
The method uses a special exponential solution of the Laplace
equation.  The solution takes the form
$e^{-\tau\,t}e^{\tau x\cdot(\omega+i\,\omega^{\perp})}$
where $\tau(>0)$ and $t$ are parameters; both $\omega$ and $\omega^{\perp
}$ are unit vectors
and satisfy $\omega\cdot\omega^{\perp}=0$.
The solution divides the space into two half planes which
have a line $\{x\,\vert\,x\cdot\omega=t\}$ as the common boundary. In one part
$\{x\,\vert\,x\cdot\omega>t\}$ the solution is growing as
$\tau\longrightarrow\infty$ and in another part $\{x\,\vert\,x\cdot\omega<t\}$ decaying.
Using this solution, we define the so-called indicator function $I_{\omega,\omega^{\perp}}(\tau,t)$ of the independent
variable $\tau$ with parameter $t$:
$$\displaystyle
I_{\omega,\omega^{\perp}}(\tau,t)
=e^{-\tau \,t}\int_{\partial\Omega}\left\{-\frac{\partial}{\partial\nu}e^{\tau x\cdot(\omega+i\omega^{\perp})
}u+\frac{\partial u}{\partial\nu}
e^{\tau x\cdot(\omega+i\omega^{\perp})
}\right\}ds.
$$
The enclosure method gives us information about the position of half plane $x\cdot\omega>t$ relative to $D$
by checking the asymptotic behaviour
of the indicator function as $\tau\longrightarrow\infty$.
For the description of the behaviour we recall
the support function $h_D(\omega)=\sup_{x\in\,D}x\cdot\omega$.
Moreover we say that $\omega$ is {\it regular} if the set
$\{x\,\vert\,x\cdot\omega=h_D(\omega)\}\cap\partial D$
consists of only one point.

What we established in \cite{I1} is: for regular $\omega$
there exist positive constants $A$ and $\mu(>1/2)$ such that,
as $\tau\longrightarrow\infty$
$$\displaystyle
\vert I_{\omega,\omega^{\perp}}(\tau,0)\vert\sim
\frac{A}{\tau^{\mu}}e^{\tau h_D(\omega)}
\tag {1.2}
$$
provided
$$\displaystyle
\text{diam}\,D<\text{dis}\,(D,\partial\Omega).
\tag {1.3}
$$
This fact is the core of the enclosure method.
Since we have the trivial identity
$$\displaystyle
I_{\omega,\omega^{\perp}}(\tau,t)
=e^{-\tau t}
I_{\omega,\omega^{\perp}}(\tau,0),
$$
from (1.2) one could conclude that:
if $t>h_D(\omega)$, then the indicator function is
decaying
exponentially; if $t=h_D(\omega)$, then the indicator function is decaying truly algebraically;
if $t<h_D(\omega)$, then the indicator function is growing exponentially.
Moreover from (1.2), we immediately obtain also the {\it one line} formula
$$\displaystyle
\lim_{\tau\longrightarrow\infty}\frac{\log\vert I_{\omega,\omega^{\perp}}(\tau,0)\vert}{\tau}
=h_D(\omega).
$$

\noindent However this is the case when the background
conductivity is {\it known}.

Consider the case when the background conductivity is given by a positive constant $\gamma$. In
this case the indicator function should be replaced with
$$\displaystyle
I_{\omega,\omega^{\perp}}(\tau,t)
=e^{-\tau \,t}\int_{\partial\Omega}\left\{-\gamma\frac{\partial}{\partial\nu}e^{\tau x\cdot(\omega+i\omega^{\perp})
}u+\gamma\frac{\partial u}{\partial\nu}
e^{\tau x\cdot(\omega+i\omega^{\perp})
}\right\}ds.
$$
Needless to say we obtain the same result as above if $\gamma$ is {\it known}.  However, if $\gamma$ is {\it unknown}, then the term
$$\displaystyle
e^{-\tau \,t}\int_{\partial\Omega}\gamma\frac{\partial}{\partial\nu}e^{\tau x\cdot(\omega+i\omega^{\perp})}uds
$$
becomes unknown and therefore one can use only the term
$$\displaystyle
e^{-\tau \,t}\int_{\partial\Omega}\gamma\frac{\partial u}{\partial\nu}
e^{\tau x\cdot(\omega+i\omega^{\perp})}
ds
\tag {1.4}
$$
if $u=f$ on $\partial\Omega$ is given.

The purpose of this paper is to give a remark on the problem: can
one still extract information about the location and shape of $D$
from the quantity (1.4) in the case when $f$ is given?

In this paper we show that, in two dimensions if the domain
occupied by the background body is enclosed by an {\it ellipse},
then it is still possible to extract some information about the location of
unknown cavities or inclusions embedded in the body without
knowing the background conductivity provided the Fourier series
expansion of the voltage on the boundary does not
contain high frequency parts (band limited) and satisfies a non
vanishing condition of a quantity involving the Fourier
coefficients.

\section{Extraction formulae}

Let $\Omega$ be the domain enclosed by an ellipse.
By choosing a suitable system of orthogonal coordinates one can write
$$\displaystyle
\Omega=\{(x_1,x_2)\,\vert\,\left(\frac{x_1}{a}\right)^2
+\left(\frac{x_2}{b}\right)^2<1\}
$$
where $a\ge b>0$.  In what follows we always use this coordinates system.

Given $\omega=(\omega_1,\omega_2)\in S^1$ set $\omega^{\perp}=(\omega_2,-\omega_1)$.
Then $x\cdot(\omega+i\omega^{\perp})=(x_1-ix_2)(\omega_1+i\omega_2)$.
Let $v=e^{\tau x\cdot(\omega+i\omega^{\perp})}$.

\subsection{Preliminary computation}

In this subsection first given $f=u\vert_{\partial\Omega}$ we
study the asymptotic behaviour of the integral
$$\displaystyle
\int_{\partial\Omega}\gamma\frac{\partial u}{\partial\nu}vds.
$$
However, integration by parts yields
$$\displaystyle
\int_{\partial\Omega}\gamma\frac{\partial u}{\partial\nu}vds
=\gamma\int_{\partial\Omega}u\frac{\partial v}{\partial\nu}ds
-\gamma\int_{\partial D}u\frac{\partial v}{\partial\nu}ds
\tag {2.1}
$$
and we have already studied the asymptotic behaviour of the second term
as described in Introduction (see (1.2)).  Therefore it suffices to study that of the first term.
Since
$$\displaystyle
\int_{\partial\Omega} u\,\frac{\partial v}{\partial\nu}ds
=\tau(\omega_1+i\omega_2)\int_{\partial\Omega}u\,v\,(\nu_1-i\nu_2)ds,
\tag {2.2}
$$
we compute the integral in the right hand side.

Write
$$\displaystyle
f(\theta)=f(a\cos\theta,b\sin\theta)
=\frac{1}{2}\alpha_0+\sum_{m=1}^{\infty}(\alpha_m\cos m\theta+\beta_m\sin m\theta)
$$
where
$$\displaystyle
\alpha_m=\frac{1}{\pi}\int_0^{2\pi}f(a\cos\theta,b\sin\theta)\cos m\theta d\theta,\,
\beta_m=\frac{1}{\pi}\int_0^{2\pi}f(a\cos\theta,b\sin\theta)\sin m\theta d\theta.
$$
Define
$$\displaystyle
\gamma_0=\alpha_0/2,\,\gamma_m=(\alpha_m-i\beta_m)/2, \gamma_{-m}=\overline{\gamma_m}, m\ge 1.
$$

\proclaim{\noindent Lemma 2.1.}
We have:
if $a=b$, then
$$\displaystyle
\int_{\partial\Omega}u\,v\,(\nu_1-i\nu_2)ds
=2\pi a^2\sum_{m=0}^{\infty}\frac{\{a\tau(\omega_1+i\omega_2)\}^m}{m!}\gamma_{m+1};
\tag {2.3}
$$
if $a>b$, then
$$
\displaystyle
\int_{\partial\Omega}u\,v\,(\nu_1-i\nu_2)ds
=2\pi ab\sum_{m=0}^{\infty}i^mJ_m(-i\sqrt{a^2-b^2}\tau(\omega_1+i\omega_2))
C_m(f)
\tag {2.4}
$$
where $C_0(f)=A_{-}\overline{\gamma_{1}}+A_{+}\gamma_1$, for $m=1,2,\cdots$
$$\displaystyle
C_m(f)=
(A_{-}\gamma_{m-1}+A_{+}\gamma_{m+1})
\left(\sqrt{\frac{a+b}{a-b}}\right)^m
+(A_{-}\overline{\gamma_{m+1}}+A_{+}\overline{\gamma_{m-1}})
\left(\sqrt{\frac{a-b}{a+b}}\right)^m
$$
and
$$\displaystyle
A_{\pm}=\frac{1}{2}\left(\frac{1}{a}\pm\frac{1}{b}\right).
$$
\endproclaim

{\it\noindent Proof.}
Set $z=e^{i\theta}$.
Since
$$\displaystyle
\nu(a\cos\theta,b\sin\theta)=\frac{1}{\displaystyle\sqrt{\left(\frac{\cos\theta}{a}\right)^2+\left(\frac{\sin\theta}{b}\right)^2}}
\left(\frac{\cos\theta}{a},\frac{\sin\theta}{b}\right)
$$
and
$$\displaystyle
ds=ab\sqrt{\left(\frac{\cos\theta}{a}\right)^2+\left(\frac{\sin\theta}{b}\right)^2}d\theta,
$$
we have
$$\displaystyle
(\nu_1-i\nu_2)ds
=ab(A_{-}z+A_{+}z^{-1})\frac{dz}{iz}.
$$
Note also that
$$\displaystyle
f(a\cos\theta,b\sin\theta)
=\sum_{m}\gamma_m z^m
$$
and
$$\displaystyle
x_1-ix_2=B_{-}z+B_{+}z^{-1}
$$
where
$$
\displaystyle
B_{\pm}=\frac{a\pm b}{2}.
$$
Using those expressions, we can write
$$\begin{array}{c}
\displaystyle
\int_{\partial\Omega}u\,v\,(\nu_1-i\nu_2)ds\\
\\
\displaystyle
=\frac{ab}{i}\sum_{m}\gamma_m
\int_{\vert z\vert=1}
(A_{-}z+A_{+}z^{-1})z^{m-1}\exp\left\{\tau(B_{-}z+B_{+}z^{-1})(\omega_1+i\omega_2)\right\}dz.
\end{array}
$$

Define
$$\displaystyle
I_l(\tau)
=\int_{\vert z\vert=1}
z^l
\exp\left\{\tau(B_{-}z+B_{+}z^{-1})(\omega_1+i\omega_2)\right\}dz.
$$

Consider the case when $a>b$.
Using the generating function of the Bessel functions, we have
$$\displaystyle
\exp\left\{\tau(B_{-}z+B_{+}z^{-1})(\omega_1+i\omega_2)\right\}
=
\sum_{n}J_n
\left(-i\sqrt{a^2-b^2}\tau(\omega_1+i\omega_2)\right)
\left(i\sqrt{\frac{a-b}{a+b}}\right)^nz^n
$$
and therefore
$$\displaystyle
I_l(\tau)
=
2\pi i(-1)^{l+1}J_{l+1}\left(-i\sqrt{a^2-b^2}\tau(\omega_1+i\omega_2)\right)
\left(-i\sqrt{\frac{a+b}{a-b}}\right)^{l+1}.
$$

\noindent
If $a=b$, then
$$\displaystyle
I_l(\tau)=0,\,\,l\le -2;\,\, I_l(\tau)=2\pi i\frac{\{a\tau(\omega_1+i\omega_2)\}^{l+1}}{(l+1)!},\,\,l\ge -1.
$$

Since
$$\displaystyle
\int_{\partial\Omega}u\,v\,(\nu_1-i\nu_2)ds
=\frac{ab}{i}
\sum_{m}\gamma_m\left(A_{-}I_m(\tau)+A_{+}I_{m-2}(\tau)\right),
$$
we obtain the desired conclusion.
\noindent
$\Box$

\subsection{Main result}

We denote by $E(\Omega)$ the set of all points on the segment that
connects the focal points $(-\sqrt{a^2-b^2},0)$ and
$(\sqrt{a^2-b^2},0)$ of $\Omega$.  It is easy to see that the
support function of the set $E(\Omega)$ is given by the formula
$h_{E(\Omega)}(\omega)=\sqrt{a^2-b^2}\vert\omega_1\vert$.

We say that a function $f(\theta)=f(a\cos\theta,b\sin\theta)$ of $\theta$
is {\it band limited} if there exists a natural number $N\ge 1$ such that,
for all $m\ge N+1$ the $m$-th Fourier coefficients $\alpha_m$ and $\beta_m$ of
the function vanish.  Then we know that $C_m(f)=0$ for all $m\ge N+2$.

Now we state the main result of this paper.

\proclaim{\noindent Theorem 2.1.}
Let $\gamma$ be a positive constant.
Assume that (1.3) is satisfied.
Let $\omega$ be regular with respect to $D$.
Let $f$ be band limited and $u$ be the solution of (1.1) with $u=f$ on $\partial\Omega$.

\noindent
(1) Let $a>b$.
Let $\omega$ satisfy $\omega_1\not=0$.
Let $f$ satisfy
$$\displaystyle
\sum_{m=1}^{\infty}(\text{sgn}\,\,\omega_1)^mm^2C_m(f)\not=0.
\tag {2.5}
$$
The formula
$$\displaystyle
\lim_{\tau\longrightarrow\infty}
\frac{1}{\tau}
\log\vert\int_{\partial\Omega}\gamma\frac{\partial u}{\partial\nu}vds\vert
=\max\,(h_D(\omega), h_{E(\Omega)}(\omega)),
\tag {2.6}
$$
is valid.

\noindent
(2) Let $a=b$.  Let $f$ satisfy: for some $N\ge 1$ $\alpha_m=\beta_m=0$ for all $m$ with $m\ge N+1$
and $\alpha_N^2+\beta_N^2\not=0$.  The formula
$$\displaystyle
\lim_{\tau\longrightarrow\infty}
\frac{1}{\tau}
\log\vert\int_{\partial\Omega}\gamma\frac{\partial u}{\partial\nu}vds\vert
=\max\,(h_D(\omega), 0),
\tag{2.7}
$$
is valid.

\endproclaim

$\bullet$  We say that a $D$ is {\it behind} the line $x\cdot\omega=t$ from the direction $\omega$
if the $D$ is contained in the half plane $x\cdot\omega<t$.
One important consequence of the formula (2.6) is:  one can know whether the unknown cavity $D$ is behind
the line $x\cdot\omega=h_{E(\Omega)}(\omega)$ from the direction $\omega$,
however, in that case one cannot know the line $x\cdot\omega=h_D(\omega)$ itself from the formula.
This shows the limit to extract the whole convex hull of $D$ without an additional assumption.

$\bullet$  The assumption that $f$ is band limited is just for a simplicity of the computation
and can be relaxed.  It is possible to apply directly the {\it saddle point method} to
study the asymptotic behaviour of the integrals in Lemma 2.1 for a $f$ that is not band limited.
Moreover we want to point out that in a practical situation, one cannot produce highly oscillatory
voltages on the boundary.  This is due to the limit of
numbers of electrodes attached on the boundary of the body.

$\bullet$  A typical example of a band-limited $f$ that satisfies (2.5) for all $\omega$ with $\omega_1\not=0$
is the $f$ given by
$$\displaystyle
f(\theta)=A\cos N\theta+B\sin N\theta
$$
where $N\ge 1$ and $A^2+B^2\not=0$.  See Remark 2.1 below for this explanation.  In general we have to
choose two $f$s corresponding to whether $\omega_1>0$ or $\omega_1<0$.

{\it\noindent Proof of Theorem 2.1.}  When $a=b$, the (2.7) is an easy consequence of (1.2),
(2.1), (2.2) and (2.3). The problem is the case when $a>b$. We
employ the {\it compound asymptotic expansion} (see page 118 of
\cite{O} for the notion of the compound asymptotic expansion) of
the Bessel function due to Hankel(see (9.09) and 9.3 of page 133
in \cite{O}):
$$\begin{array}{l}
\displaystyle
J_{m}(z)
\sim
\left(\frac{2}{\pi z}\right)^{1/2}\times
\\
\displaystyle
\left\{
\cos\left(z-\frac{m\pi}{2}-\frac{\pi}{4}\right)
\sum_{s=0}^{\infty}(-1)^s\frac{A_{2s}(m)}{z^{2s}}
-\sin\left(z-\frac{m\pi}{2}-\frac{\pi}{4}\right)
\sum_{s=0}^{\infty}(-1)^s\frac{A_{2s+1}(m)}{z^{2s+1}}\right\}
\end{array}
\tag {2.8}
$$
as $z\longrightarrow\infty$ in $\vert\text{arg}\,z\vert\le\pi-\delta$
for each fixed $\delta\in\,]0,\,\pi[$ where $A_0(m)=1$ and,
for $s=1,2,\cdots$
$$\displaystyle
A_s(m)=\frac{1}{s!8^s}
(4m^2-1^2)(4m^2-3^2)\cdots(4m^2-(2s-1)^2).
$$
First we consider the case when $\omega_1>0$.  From (2.8) in the case when
$z=-i\sqrt{a^2-b^2}\tau(\omega_1+i\omega_2)$ we obtain
$$\displaystyle
J_m(z)=\left(\frac{1}{2\pi z}\right)^{1/2}
e^{iz}(-i)^me^{-i\pi/4}
\left(1-\frac{4m^2-1}{8iz}+O(\frac{1}{\tau^2})\right)
\tag {2.9}
$$
as $\tau\longrightarrow\infty$.
Since $f$ is band limited, one can find $N\ge 1$
such that, for all $m\ge N+1$ the $m$-th Fourier coefficients $\alpha_m$ and $\beta_m$ of $f$ vanish.
Then $C_m(f)=0$ for $m\ge N+2$ and from (2.4) and (2.9) we obtain
$$\begin{array}{c}
\displaystyle
\int_{\partial\Omega}uv(\nu_1-i\nu_2)ds
=2\pi ab
\left(\frac{1}{2\pi z}\right)^{1/2}
e^{iz}e^{-i\pi/4}\\
\\
\displaystyle
\times
\left\{\left(1+\frac{1}{8iz}\right)\sum_{m=0}^{N+1}C_m(f)
+i\frac{1}{2z}\sum_{m=1}^{N+1}m^2 C_m(f)+O(\frac{1}{\tau^2})\right\}.
\end{array}
\tag {2.10}
$$
Here we claim that
$$\displaystyle
\sum_{m=0}^{N+1}C_m(f)=0.
\tag {2.11}
$$
It suffices to prove the claim in the case when
$$\displaystyle
f(a\cos\theta, b\sin\theta)=\alpha_j\cos j\theta+\beta_j\sin j\theta
\tag {2.12}
$$
for each fixed $j=1,2,\cdots, N$.  Since $\sum_{m=0}^{\infty}C_m(f)
=C_{j-1}(f)+C_j(f)+C_{j+1}(f)$ and we have
$$\begin{array}{c}
\displaystyle
C_{j+1}(f)=A_{-}\gamma_j\left(\sqrt{\frac{a+b}{a-b}}\,\right)^{j+1}
+A_{+}\overline{\gamma_j}\left(\sqrt{\frac{a-b}{a+b}}\,\right)^{j+1},\\
\\
C_j(f)=0,\\
\\
\displaystyle
C_{j-1}(f)=A_{+}\gamma_j\left(\sqrt{\frac{a+b}{a-b}}\,\right)^{j-1}
+A_{-}\overline{\gamma_j}\left(\sqrt{\frac{a-b}{a+b}}\,\right)^{j-1},
\end{array}
$$
we get
$$
\displaystyle
\sum_{m=0}^{\infty}C_m(f)
=\left\{A_{+}+A_{-}\left(\frac{a+b}{a-b}\right)\right\}
\left\{\gamma_j\left(\sqrt{\frac{a+b}{a-b}}\,\right)^{j-1}
+\overline{\gamma_j}\left(\sqrt{\frac{a-b}{a+b}}\,\right)^{j+1}\right\}.
$$
Since
$$\displaystyle
A_{+}+A_{-}\left(\frac{a+b}{a-b}\right)=0,
$$
we see that the claim (2.11) is valid.
Therefore (2.10) becomes
$$
\displaystyle
\int_{\partial\Omega}uv(\nu_1-i\nu_2)ds
=i\pi abz^{-1}
\left(\frac{1}{2\pi z}\right)^{1/2}
e^{iz}e^{-i\pi/4}
\left(\sum_{m=1}^{N+1}m^2 C_m(f)+O(\frac{1}{\tau})\right).
\tag {2.13}
$$
Set $\omega_1+i\omega_2=e^{i\vartheta}$ with
$-\pi/2<\vartheta<\pi/2$. Then
$z^{1/2}=\sqrt{\tau}(a^2-b^2)^{1/4}e^{i(\vartheta-\pi/2)/2}$.
Since $e^{iz}=e^{\tau
h_{E(\omega)}(\omega)}e^{i\tau\sqrt{a^2-b^2}\omega_2}$, from
(1.2), (2.1), (2.2) and (2.13) we obtain the compound asymptotic
formula:
$$\begin{array}{c}
\displaystyle
\int_{\partial\Omega}\gamma\frac{\partial u}{\partial\nu}vds\\
\\
\displaystyle
\sim-\gamma
\sqrt{\frac{\pi}{2}}ab(a^2-b^2)^{-3/4}e^{-i\vartheta/2}\tau^{-1/2}
e^{\tau h_{E(\Omega)}(\omega)}e^{i\tau\sqrt{a^2-b^2}\omega_2}
\sum_{m=1}^{N+1}m^2 C_m(f)
-\gamma e^{\tau h_D(\omega)}\frac{A}{\tau^{\mu}}.
\end{array}
$$
From this we know that the quantity
$$\displaystyle
\exp\left\{-\tau\max\,(h_D(\omega),h_{E(\Omega)}(\omega))\right\}\vert\int_{\partial\Omega}
\gamma\frac{\partial u}{\partial\nu}vds\vert
$$
is truly {\it algebraic} decaying as $\tau\longrightarrow\infty$.  
Note that we have used
the lower bound of $\mu$: $\mu>1/2$.
Therefore we obtain the formula (2.6).
Next consider the case when $\omega_1<0$.  Write
$\displaystyle
R_{\omega}(\tau;f)=\int_{\partial\Omega}fv(\nu_1-i\nu_2)ds$. Then
we have $R_{\omega}(\tau;f)=-R_{-\omega}(\tau;f^*)$ where
$f^*(x)=f(-x)$.  Since the $m$-th Fourier coefficients of $f^*$ are given
by $(-1)^m$ times those of $f$ and the first
component of $-\omega$ is positive, we can derive the
corresponding result in the case when $\omega_1<0$ from the result
in the case when $\omega_1>0$ by replacing $C_m(f)$ in the condition
(2.5) with $-(-1)^mC_m(f)$. $\Box$

{\bf\noindent Remark 2.1.}
Fix $j=1,2,\cdots,N$ and let $f$ be given by (2.12).
Then a direct computation similar to the proof of the claim (2.11) yields
$$\begin{array}{c}
\displaystyle
\sum_{m=1}^{\infty}m^2C_m(f)
=(j-1)^2C_{j-1}(f)+(j+1)^2C_{j+1}(f)\\
\\
\displaystyle
=-\frac{2}{ab}j(a^2-b^2)^{-(j-1)/2}
\{(a+b)^j\gamma_j-(a-b)^j\overline{\gamma_j}\}.
\end{array}
$$
This yields also
$$\begin{array}{c}
\displaystyle
\sum_{m=1}^{\infty}(-1)^mm^2C_m(f)
=(-1)^{j-1}\sum_{m=1}^{\infty}m^2C_m(f)\\
\\
\displaystyle
=(-1)^{j}
\frac{2}{ab}j(a^2-b^2)^{-(j-1)/2}
\{(a+b)^j\gamma_j-(a-b)^j\overline{\gamma_j}\}.
\end{array}
$$
These yield: a $f$ whose Fourier coefficients $\alpha_j$ and $\beta_j$
vanish for all $j\ge N+1$ with some $N\ge 1$, satisfies the condition (2.5)
if and only if
$$\displaystyle
\sum_{j=1}^{N} (\text{sgn}\,\omega_1)^jj(a^2-b^2)^{-(j-1)/2}
\{(a+b)^j\gamma_j-(a-b)^j\overline{\gamma_j}\}\not=0.
\tag {2.14}
$$
It is clear that there are many $f$s satisfying the condition (2.14).

{\bf\noindent Remark 2.2.}
In (1) the case when $\omega_1=0$ is not treated. In this case $\omega_2=\pm 1$.
If $\omega_2=1$, then from (2.8) we have
$$\begin{array}{l}
\displaystyle
J_m(z)=\left(\frac{1}{2\pi\sqrt{a^2-b^2}\tau}\right)^{1/2}
\times
\\
\\
\displaystyle
\left\{e^{i\tau\sqrt{a^2-b^2}}(-i)^me^{-i\pi/4}
\left(1+i\frac{4m^2-1}{8z}\right)
+e^{-i\tau\sqrt{a^2-b^2}}i^me^{i\pi/4}
\left(1-i\frac{4m^2-1}{8z}\right)
\right\}+O(\tau^{-5/2})
\end{array}
$$
where $z=-i\sqrt{a^2-b^2}\,\tau(\omega_1+i\omega_2)$.  Then from (1.2), (2.1), (2.2) and (2.4) the problem can be reduced
to the study of the asymptotic behaviour of the quantity
$$\displaystyle
\sum_{m=0}^{N+1}
\left\{e^{i\tau\sqrt{a^2-b^2}}e^{-i\pi/4}
\left(1+i\frac{4m^2-1}{8z}\right)
+(-1)^m e^{-i\tau\sqrt{a^2-b^2}}e^{i\pi/4}
\left(1-i\frac{4m^2-1}{8z}\right)
\right\}C_m(f)
\tag {2.15}
$$
as $\tau\longrightarrow\infty$.  This seems very complicated for general $\tau$.  However,
if we choose
$$\displaystyle
\tau=\frac{l\pi}{\sqrt{a^2-b^2}},\,\,l=1,2,\cdots,
\tag {2.16}
$$
then (2.15) becomes
$$\begin{array}{c}
\displaystyle
(-1)^le^{-i\pi/4}\left\{\sum_{m=0}^{N+1}\left(1+i\frac{4m^2-1}{8z}\right)
C_m(f)+i\sum_{m=0}^{N+1}(-1)^m
\left(1-i\frac{4m^2-1}{8z}\right)
C_m(f)\right\}\\
\\
\displaystyle
=\frac{(-1)^le^{-i\pi/4}i}{2z}\sum_{m=1}^{N+1}m^2\{C_m(f)-iC_m(f^*)\}.
\end{array}
$$
Note that we have used the claim (2.11) for $f$ and $f^*$.
Therefore if $f$ satisfies the condition
$$\displaystyle
\sum_{m=1}^{\infty}m^2\{C_m(f)-iC_m(f^*)\}\not=0
\tag {2.17}
$$
instead of (2.5), then for $\tau$ given by (2.16), the formula
$$\displaystyle
\lim_{l\longrightarrow\infty}
\frac{1}{\tau}
\log\vert\int_{\partial\Omega}\gamma\frac{\partial u}{\partial\nu}vds\vert
=\max\,(h_D(\omega), 0),
$$
is valid.   By replacing $f$ with $f^*$, we know also that: if $\omega_2=-1$, then the same formula is valid
provided
$$\displaystyle
\sum_{m=1}^{\infty}m^2\{C_m(f)+iC_m(f^*)\}\not=0
\tag {2.18}
$$
instead of (2.17).  From the computation in Remark 2.1 one can sum
the conditions (2.17) and (2.18) up in the single form:
$$\displaystyle
\sum_{j=1}^{N}
\left\{1+(-1)^{j}(\text{sgn}\,\omega_2)i\right\}
j(a^2-b^2)^{-(j-1)/2}
\{(a+b)^j\gamma_j-(a-b)^j\overline{\gamma_j}\}\not=0
$$
where $N\ge 1$ and chosen in such a way that, for all $m\ge N+1$ the $m$-th Fourier coefficients of $f$ vanish.

\subsection{Uniqueness}

As a corollary of Theorem 2.1 we obtain a uniqueness theorem.

\proclaim{\noindent Corollary 2.1.} Let $\gamma$ be a positive
constant. Assume that $D$ satisfies (1.3).

(1)  Let $\Omega$ be a domain enclosed by an ellipse.  Let $f_{+}$ and $f_{-}$ be band limited
and satisfy
$$\displaystyle
\sum_{m=1}^{\infty}(\pm)^m m^2 C_m(f_{\pm})\not=0.
$$
Let $u_{\pm}$ be the solution of (1.1) with $u_{\pm}=f_{\pm}$ on $\partial\Omega$.
Then the Neumann data $\gamma\partial u_{+}/\partial\nu$
and $\gamma\partial u_{-}/\partial\nu$ on $\partial\Omega$
uniquely determine the convex hull of $D\cup E(\Omega)$.

(2) Let $\Omega$ be a domain enclosed by a circle.  Let $f$ be
band limited and non constant.  Let $u$ be the solution of (1.1)
with $u=f$ on $\partial\Omega$. Then the Neumann data
$\gamma\partial u/\partial\nu$ uniquely determines the convex hull
of $D\cup\{0\}$.

\endproclaim

\noindent We emphasize that $\gamma$ is {\it unknown}. This makes
the situation difficult definitely. Assume that we have two
unknowns $(D,\gamma)=(D_1, \gamma_1), (D_2,\gamma_2)$ and
solutions $u_1$ and $u_2$ both satisfying (1.1) and the boundary
condition $u=f$ on $\partial\Omega$. The key point of a standard
and traditional approach is to prove that if $\gamma_1\partial
u_1/\partial\nu=\gamma_2\partial u_2/\partial\nu$ on
$\partial\Omega$, then $u_1=u_2$ in a neighbourhood of
$\partial\Omega$. If $\gamma_1=\gamma_2$, then the conclusion is
true because of the uniqueness of the Cauchy problem for the
Laplace equation. However, if $\gamma$ is unknown, i.e., the
assumption $\gamma_1=\gamma_2$ is dropped, one can not immediately
get the conclusion (note that we are considering a finite set of
observation data not the full Dirichlet-to-Neumann map). Our
approach skips this point by using an analytical formula that
directly connects the data with unknown discontinuity.

The proof of Corollary 2.1 is based on:  given $D$ the set of all
directions that are not regular with respect to $D$ is a finite
set; the formulae (2.6) are valid for $f=f_{\pm}$ in (1); the
formula (2.7) is valid for $f$ in (2). Therefore, for example, in
(1) we see that the Neumann data uniquely determine the values of
$\max\,(h_D(\omega), h_{E(\Omega)}(\omega))$ which is the support
function of the convex hull of $D\cup E(\Omega)$ at the directions
$\omega$ except for a finite set of directions. Since the support
function $h_D$ and $h_{E(\Omega)}$ are continues on the unit
circle and so is
$\max\,(h_D(\,\cdot\,),h_{E(\Omega)}(\,\cdot\,))$.  A density
argument yields the desired uniqueness.

{\bf\noindent Remark 2.3.}  If $\partial D$ is {\it smooth},  then
(2) of Corollary 2.1 does not hold. Let $\Omega$ be the unit open
disc centered at the origin of the coordinates system and for
$0<R<1$ let $D(R)$ be the open disc centered at the origin with
the radius $R$. Let $0<R_1, R_2<1$.  Fix an integer $m\ge 1$. For
each $j=1,2$ let $u_j$ be the weak solution of  the problem (1.1)
with $D=D(R_j)$ and the Dirichlet data
$u_j(r,\theta)\vert_{r=1}=\cos\,m\theta$ where $(r,\theta)$
denotes the usual polar coordinates centered at the origin. Then
we know that
$$\displaystyle
u_1(r,\theta)=\frac{1}{1+R_1^{2m}}(r^m+R_1^{2m}r^{-m})\cos\,m\theta,\,\,
u_2(r,\theta)=\frac{1}{1+R_2^{2m}}(r^m+R_2^{2m}r^{-m})\cos\,m\theta.
$$
This yields
$$\displaystyle
\frac{1+R_2^{2m}}{1-R_2^{2m}}\frac{\partial u_2}{\partial\nu}=m\cos\,m\theta
=\frac{1+R_1^{2m}}{1-R_1^{2m}}\frac{\partial u_1}{\partial\nu}\,\,\text{on}\,\partial\Omega.
$$
Since $R_1$ and $R_2$ are arbitrary chosen, this means that one cannot uniquely determine
$D(R)$ from the single set of the Dirichlet and Neumann data $f(\theta)=\cos\,m\theta$ and $\gamma\partial u/\partial\nu$
on $\partial\Omega$ in the case when $\gamma=(1+R^{2m})(1-R^{2m})$.  This suggests that the {\it singularity} of $\partial D$
is essential for the validity of (2) in Corollary 2.1.

\section{An application to the inverse conductivity problem}

The idea in the proof of Theorem 2.1 can be applied to the case
when the unknown domain $D$ is a model of an {\it inclusion}.

We assume that the conductivity $k=k(x)$ of the body that occupies
$\Omega$ is given by $k(x)=\gamma$ if $x\in\Omega\setminus D$;
$k(x)=\tilde{\gamma}$ if $x\in D$.  It is assumed that the
$\gamma$ and $\tilde{\gamma}$ are positive constants and satisfy
$\gamma\not=\tilde{\gamma}$.  The voltage potential $u$ inside the
body satisfies the equation $\displaystyle\nabla\cdot k\nabla u=0$
in $\Omega$.
Given $\omega=(\omega_1,\omega_2)\in S^1$ set $\omega^{\perp}=(\omega_2,-\omega_1)$.
Let $\tau>0$ and $v=e^{\tau x\cdot(\omega+i\omega^{\perp})}$.

In \cite{Ie2} we have already proved that if $u$ is
not a constant function and $D$ is polygonal and satisfies the
condition (1.3), then for a given direction $\omega$ that is
regular with respect to $D$ the formula
$$\displaystyle
\lim_{\tau\longrightarrow\infty}
\frac{1}{\tau}\log\vert\int_{\partial\Omega}\left(\gamma\frac{\partial u}{\partial\nu}v-
\gamma\frac{\partial v}{\partial\nu}u\right)ds\vert
=h_D(\omega),
$$
is valid. Note that $k=\gamma$ on $\partial\Omega$ and we do not
assume that the conductivity $\tilde{\gamma}$ of $D$ is known.

Here we propose the same question as that of Introduction.
Assume that
we do not know $k$ in the whole domain. Given a non constant
voltage potential $f=u\vert_{\partial\Omega}$ on $\partial\Omega$
is it possible to extract some information about the location of
$D$ from the corresponding current density $k\partial
u/\partial\nu$ on $\partial\Omega$?

The answer is yes in the case
when the $\Omega$ is enclosed by an ellipse.
It starts with recalling the equation
$$\displaystyle
\int_{\partial\Omega} \gamma\frac{\partial u}{\partial\nu} vds
=\int_{\partial\Omega}\gamma\frac{\partial v}{\partial\nu}uds-(\gamma-\tilde{\gamma})
\int_{\partial D}u\frac{\partial v}{\partial\nu}ds.
\tag {3.1}
$$
Recall Key Lemma in \cite{Ie2}:  there exist positive constants $B$ and $\lambda(>1/2)$ such that,
as $\tau\longrightarrow\infty$
$$\displaystyle
\vert \int_{\partial D}u\frac{\partial v}{\partial\nu}ds\vert\sim
\frac{B}{\tau^{\lambda}}e^{\tau h_D(\omega)}.
\tag {3.2}
$$
Then from (2.2), (3.1), (3.2) and Lemma 2.1 we see that
the completely same statements as those in Theorem 2.1, Corollary 2.1 and Remarks 2.1 and 2.2 are valid.

{\bf\noindent Remark 3.1.}  In \cite{Ie3} we employed the {\it difference} of the values of the voltage
at arbitrary fixed two points on the boundary of a {\it general} two-dimensional bounded domain $\Omega$ with
smooth boundary.  More precisely we introduced the operator
$$\displaystyle
\Lambda_{k}(P,Q):g\longmapsto u(P)-u(Q)
$$
where $P$ and $Q$ are two arbitrary points on $\partial\Omega$;
$g$ satisfies $\displaystyle\int_{\partial\Omega}gds=0$; the $u$ is a solution
of the equation $\nabla\cdot k\nabla u=0$ in $\Omega$ and
satisfies the Neumann boundary condition $k\partial
u/\partial\nu=g$ on $\partial\Omega$.

Given $\omega=(\omega_1,\omega_2)\in S^1$ set $\omega^{\perp}=(\omega_2,-\omega_1)$.
Let $\tau>0$ and $v=e^{\tau x\cdot(\omega+i\omega^{\perp})}$.
What we have proved is:  if $g=\partial v/\partial\nu$ on $\partial\Omega$ and
$D$ is polygonal and satisfies the
condition (1.3), then for a given direction $\omega$ that is
regular with respect to $D$ the formula
$$\displaystyle
\lim_{\tau\longrightarrow\infty}
\frac{1}{\tau}\log\vert\left\{\Lambda_{k}(P,Q)-\Lambda_{\gamma}(P,Q)\right\}(g)\vert
=h_D(\omega),
\tag {3.3}
$$
is valid.  Note that we have used the relationship
$$
\displaystyle
\left\{\Lambda_{k}(P,Q)-\Lambda_{\gamma}(P,Q)\right\}(g)
=\frac{1}{\gamma}\left\{\Lambda_{k/\gamma}(P,Q)-\Lambda_1(P,Q)\right\}(g).
$$

If $\gamma$ is unknown, then one cannot use the term
$\Lambda_{\gamma}(P,Q)(g)$ in (3.3).  However, that has the simple form
$$\displaystyle
\Lambda_{\gamma}(P,Q)(g)=\frac{1}{\gamma}\left\{v(P)-v(Q)\right\}
$$
for $g=\partial v/\partial\nu$ on $\partial\Omega$.
Using this form, Proposition 3.1 and Lemma 3.1 in \cite{Ie3}, one immediately gets the following
formulae provided $D$ is polygonal and satisfies the
condition (1.3) and $\omega$ is
regular with respect to $D$:

$\bullet$  if $\omega$ is not perpendicular to the line passing through $P$ and $Q$, then
$$\displaystyle
\lim_{\tau\longrightarrow\infty}
\frac{1}{\tau}\log\vert\Lambda_{k}(P,Q)(g)\vert
=\max\,\left(h_D(\omega), h_{\{P,\,Q\}}(\omega)\right);
$$

$\bullet$  if $\omega$ is perpendicular to the line passing through $P$ and $Q$,
choose, for example,
$$\displaystyle
\tau=\frac{\pi}{\vert P-Q\vert}\left(\frac{1}{2}+2l\right),\,\,l=0,1,2,\cdots,
$$
then
$$\displaystyle
\lim_{l\longrightarrow\infty}
\frac{1}{\tau}\log\vert\Lambda_{k}(P,Q)(g)\vert
=\max\,\left(h_D(\omega), h_{\{P,\,Q\}}(\omega)\right).
$$

\section{Conclusion}

We confirmed that: in the
case when the background conductivity is homogeneous and {\it
unknown} the enclosure method still works provided:

$\bullet$  the domain that is occupied by a background
body has a simple geometry;

$\bullet$  the Fourier series
expansion of the voltage on the boundary does not
contain high frequency parts (band limited) and satisfies a non
vanishing condition of a quantity involving the Fourier
coefficients.

However, the method yields a less information about the location
and shape of unknown cavity or inclusion compared with the case
when the conductivity is {\it known}.  We found an explicit {\it
obstruction} that depends on the geometry of the background body.

$$\quad$$

\centerline{{\bf Acknowledgements}}

This research was partially supported by Grant-in-Aid for
Scientific Research (C)(No.  18540160) of Japan  Society for the
Promotion of Science.

\end{document}